\documentclass{article}

\usepackage{amsmath,amssymb,amsfonts,amsthm,latexsym}
\usepackage{hyperref, enumerate}
\usepackage{tikz}
\usepackage{rotating}
\usepackage{ifthen}
\usepackage{geometry}

\usetikzlibrary{calc, math}
\usetikzlibrary{decorations.markings}

\newtheorem{example}{\bf Example}
\newtheorem{conjecture}{\bf Conjecture}

\newtheorem{proposition}{\bf Proposition}

\newtheorem{theorem}{\bf Theorem}
\newtheorem{definition}{\bf Definition}
\newtheorem*{definition*}{\bf Definition}
\newtheorem{problem}{\bf Problem}
\newtheorem*{problem*}{\bf Problem}

\newtheorem*{comment*}{\bf Comment}

\allowdisplaybreaks

\title{3-critical subgraphs of snarks}

\author{Imran Allie}

\begin{document}

\maketitle

\begin{abstract}
	In this paper we further our understanding of the structure of class two cubic graphs, or snarks, as they are commonly known. We do this by investigating their 3-critical subgraphs, or as we will call them, minimal conflicting subgraphs. We consider how the minimal conflicting subgraphs of a snark relate to its possible minimal 4-edge-colourings. We fully characterise the relationship between the resistance of a snark and the set of minimal conflicting subgraphs. That is, we show that the resistance of a snark is equal to the minimum number of edges which can be selected from the snark, such that the selection contains at least one edge from each minimal conflicting subgraph. We similarly characterise the relationship between what we call \textit{the critical subgraph} of a snark and the set of minimal conflicting subgraphs. The critical subgraph being the set of all edges which are conflicting in some minimal colouring of the snark. Further to this, we define groups, or \textit{clusters}, of minimal conflicting subgraphs. We then highlight some interesting properties and problems relating to clusters of minimal conflicting subgraphs. 
\end{abstract}

\section{Introduction}{\label{Introduction}}
%The edge chromatic number of a graph is the minimum number of colours required to assign to each edge of the graph a colour, such that no two adjacent edges are assigned the same colour. As is well-known, the edge chromatic number of a cubic graph is either three or four. Such graphs are referred to as cubic class one and cubic class two graphs, respectively. Cubic class two graphs have long been of particular interest in graph theory, largely for the fact that many major problems in graph theory are easily solvable for graphs which are not cubic class two. Tutte's 5-flow conjecture \cite{tutte} and the cycle double cover conjecture \cite{szekeres} are major examples of the these problems. In fact, for some major problems, only a subset of cubic class two graphs present particular difficulty. That is, a subset called \textit{snarks}. Snarks are, in addition, cyclically 4-edge-connected and have girth at least 5. This restriction removes much triviality in the discussion of cubic class two graphs.

As is well-known, the edge chromatic number of a cubic graph is either three or four. Such graphs are referred to as cubic class one and cubic class two graphs, respectively. Cubic class two graphs are more commonly known as snarks. Snarks have long been of particular interest in graph theory, largely for the fact that many major problems in graph theory are easily solvable for graphs which are not snarks. Tutte’s 5-flow conjecture \cite{tutte} and the cycle double cover conjecture \cite{szekeres} are major examples of the these problems.

Let $G=(V,E)$ be a graph. A $k$-edge-colouring, $f$, of $G$ is a mapping from the set of edges of $G$ to a set of $k$ colours. That is,  $f : E \longrightarrow \{1,\dots,k\}$. $f$ is a \textit{proper} $k$-edge-colouring  of $G$ if no two adjacent elements in $E$ are mapped to the same colour. By Vizing's theorem \cite[Theorem 6.2]{bondy}, if $G$ is a graph and $f$ is a proper colouring then the smallest possible value of $k$ is $\Delta$ or $\Delta + 1$, where $\Delta$ is the maximum degree of any vertex in $G$. If the smallest possible value of $k$ is $\Delta$,  then we say that $G$ is class one, or $\Delta$-edge-colourable. Otherwise we say that $G$ is class two, or $(\Delta + 1)$-edge-colourable. Given a $k$-edge-colouring $f$, we call the set $f^{-1}(i)$ a colour class, for each $i \in \{1,\dots,k\}$. A vertex $v$ is \textit{conflicting} with regard to $f$ if more than one of the edges incident to $v$ are mapped to the same colour.

The \textit{resistance} of $G$, denoted as $r(G)$, is defined as the $\min\{|f^{-1}(i)| : f$ is a proper $(\Delta+1)$-edge-colouring of G and $f^{-1}(i)$ is a colour class$\}$. That is, the minimum number of edges that can be removed from a graph such that the resulting graph is 3-edge-colourable \cite{steffen}. As it turns out, somewhat counter-intuitively perhaps, the resistance of $G$ equals the \textit{vertex resistance} of $G$, denoted as $r_v(G)$, which is the minimum number of vertices that needs to be removed from $G$ such that the resultant graph is class one. If given a 3-edge-colouring of $G$ with $r(G)$ conflicting vertices, we can also find a proper 4-edge-colouring of $G$ with $r(G)$ edges being mapped to one particular colour. Furthermore, the conflicting vertices in the 3-edge-colouring have a one-to-one relationship with the set of edges mapped to the fourth colour in the 4-edge-colouring. This is a result we use implicitly going forward, in that we do not consider 3-edge colourings with conflicting vertices. We consider only proper 3-edge colourings or proper 4-edge-colourings of cubic graphs. It has also been proven that $r(G) = 0$ or $r(G) \geq 2$ for any cubic graph $G$ \cite{fiol}.

If $f$ is a proper $k$-edge-colouring of a graph $G$ and $|f^{-1}(i)| = r(G)$ for some $i \in \{1, \dots, k\}$, then we call $f$ a \textit{minimal colouring} \cite{steffen2}. For cubic graphs, we will use colour sets $\{1,2,3\}$ and $\{0,1,2,3\}$ for class one and class two graphs, respectively. We will assume $|f^{-1}(0)| = r(G)$ for a minimal colouring $f$ of $G$. Given a minimal colouring $f$ of $G$, if $f(e)=0$ for some edge $e \in G$ then we call $e$ a \textit{conflicting edge} with regard to $f$. Now, let $H \subset G$ and let $f_H$ be a proper colouring of $H$. A proper colouring of $G$, $f_G$, with $f_G(e) = f_H(e)$ for all $e \in H$ is called an \textit{extension} of $f_H$.%, or we say that $f_H$ can be \textit{extended} to $f_G$. 
If $f_G$ is such that the number of conflicting edges in $G-H$ is minimal given $f_H$, then we call $f_G$ a \textit{minimal extension} of $f_H$.
%, or we say that $f_H$ can be \textit{minimally extended} to $f_G$. We say that two minimal colourings of $G$ are \textit{equivalent} if there conflicting edges coincide. 

In this paper, we further understand the complexity of these graphs by extending on the definition of conflicting zones introduced in \cite{fiol}, although we opt for the term conflicting subgraphs. We define minimal conflicting subgraphs, as well as consequent concepts such as: the buffer subgraph, which is the maximum subgraph containing no edges in or adjacent to any minimal conflicting subgraph; the critical subgraph, which is the subgraph containing all edges which are conflicting in some minimal colouring of the graph; and clusters of minimal conflicting subgraphs which are essentially overlapping minimal conflicting subgraphs (formal definitions to follow). We then prove the following insight about cubic class two graphs. That for any collection of edges $R$ in a cubic graph $G$, such that $R$ contains an edge from every minimal conflicting subgraph in $G$, $G-R$ is 3-edge-colourable. Furthermore, for such an $R$ with minimal possible order, the resistance of $G$ is $|R|$. We are able to then characterise the resistance, as well as the critical subgraph of a graph, in terms of the set of minimal conflicting subgraphs. Finally, we discuss further problems of consideration.

%%%%%%%%%%%%%%%%%%%%%%%%%%%%%%%%%%%%%%%%%%%%%%%%%%%%%%%%%%%%%%
%%%%%%%%%%%%%%%%%%%%%%%%%%%%%%%%%%%%%%%%%%%%%%%%%%%%%%%%%%%%%%
%%%%%%%%%%%%%%%%%%%%%%%%%%%%%%%%%%%%%%%%%%%%%%%%%%%%%%%%%%%%%%
%%%%%%%%%%%%%%%%%%%%%%%%%%%%%%%%%%%%%%%%%%%%%%%%%%%%%%%%%%%%%%
%%%%%%%%%%%%%%%%%%%%%%%%%%%%%%%%%%%%%%%%%%%%%%%%%%%%%%%%%%%%%%
%%%%%%%%%%%%%%%%%%%%%%%%%%%%%%%%%%%%%%%%%%%%%%%%%%%%%%%%%%%%%%
\section{Minimal conflicting subgraphs}

A conflicting subgraph of a cubic graph $G$ is defined as a subgraph $H$ of $G$ which does not admit a proper 3-edge-colouring. That is essentially, a subgraph which itself is not 3-edge-colourable. This idea was introduced in \cite{fiol}, where it was called a conflicting zone. With a view to further understand what makes a cubic graph class two, we extend on this idea by defining \textit{minimal conflicting subgraphs}. The essential idea being to isolate from the graph that which is non 3-edge-colourable. 

\begin{definition} \label{mcz}{\rm
	Let $G$ be a subcubic graph and let $M$ be a conflicting subgraph of $G$. If for any $e \in E(M)$ we have that $M - {e}$ is not a conflicting subgraph, then we call $M$ a \textit{minimal conflicting subgraph} of $G$. Let 
$$M_G = \bigcup \{ M ~|~M {\rm~is~a~minimal~conflicting~subgraph~of~} G\}.$$ 
We call $M_G$ the \textit{maximal conflicting subgraph} of $G$. Let 
$$C_G = \{ e \in E(G)~|~e \notin M_G~ {\rm and}~e~{\rm is~ adjacent~ to~ some}~e'~{\rm in}~ M_G\}.$$ 
We call $C_G$ \textit{the conflict-cut set} of $G$. Let 
$$B_G = \{ e \in E(G)~|~e \notin M_G \cup C_G\}.$$ 
We call $B_G$ the \textit{buffer subgraph} of $G$.}
\end{definition}

A subcubic graph $G$ is called 3-critical if it has chromatic index 4 and $G - e$ has chromatic index 3 for every $e \in G$. It is easy to see that this definition coincides with our definition of minimal conflicting subgraphs, in that a minimal conflicting subgraph can be thought of as a 3-critical subgraph. If the 3-critical subgraphs represent only that which is essentially non 3-edge-colourable, then the buffer subgraph represents that which is essentially redundant in contributing to the non colourability of the cubic graph. We list some properties of 3-critical subgraphs, or minimal conflicting subgraphs, of subcubic graphs. First we present some known properties of 3-critical graphs in general, after which we prove some more pertinent properties for our purposes regarding minimal conflicting subgraphs. 

\begin{proposition} \label{properties_of_3crits}
	Let $M$ be a 3-critical graph. The following statements are true.
	\begin{enumerate}
		\item[(i)] 	$r(M)=1$ and every edge $e \in M$ is conflicting in some minimal colouring of $M$.
		\item[(ii)] 	$M$ is strictly subcubic.
		\item[(iii)] 	$M$ is bridgeless.
		\item[(iv)] 	Every vertex in $M$ has degree two or three.
		\item[(v)] 	Every vertex in $M$ has at least two neighbours of degree three. 
		%\item[(vi)] 	$M$ has at least three vertices of degree two. 
	\end{enumerate} 
	\begin{proof} 
		These are known properties of 3-critical graphs and we omit the proofs.
	\end{proof}
\end{proposition}

\begin{proposition} \label{properties_of_mczs}
	Let $G$ be a bridgeless cubic graph. The following statements are true.
	\begin{enumerate}
		\item[(i)] 	The distance between any two disjoint minimal conflicting subgraphs of $G$ is at least one.
		\item[(ii)]	Every conflicting subgraph in $G$ contains a minimal conflicting subgraph.
	\end{enumerate} 
	\begin{proof} \hfill

		\begin{enumerate}
			\item[(i)] This follows on directly from Proposition \ref{properties_of_3crits} (iv).
			\item[(ii)] Let $M$ be a conflicting subgraph of $G$. Choose an edge $e \in M$. We check $e$ by considering $r(M - \{ e \})$. If $r(M - \{ e \}) \neq 0$ then remove $e$ from $M$. If $r(M - \{ e \})=0$ then leave $M$ as is and mark $e$ as checked. Continue checking edges in $M$ until every edge is checked. Once every edge is checked, $M$ is then a minimal conflicting subgraph. \vspace{-1pc}
		\end{enumerate}
	\end{proof}
\end{proposition}

We begin our investigation into these structures. We consider their existence relative to conflicting edges in minimal colourings. Note that although our primary interest is in cubic graphs, some results are applicable to subcubic graphs as well and are stated as such.

\begin{proposition} \label{r_distinct_mczs}
	Let $G$ be a subcubic class two graph and let $f$ be a minimal colouring of $G$. For each conflicting edge $e$ with regard to $f$, there exists at least one minimal conflicting subgraph which contains $e$ and also contains no other conflicting edge with regard to $f$.
	\begin{proof}
		Let $f$ be a minimal colouring of $G$ and let $R = \{ e_1, \dots, e_r \}$ be the set of conflicting edges with regard to $f$. For each $i \in \{ 1, \dots ,r \}$ let $M_i = \{ e_i\}$ and conduct the following process. Choose an edge $e$ not contained in $M_i \cup R$ which is adjacent to some edge in $M_i$. Add edge $e$ to $M_i$. While $r(M_i)=0$, we keeping adding such edges. Since $r(G - (R-\{ e_i\}))$ must equal 1, we know that eventually we will have $r(M_i)=1$. If $r(M_i)=1$ then $M_i$ is a conflicting subgraph which contains no other conflicting edge with regard to $f$ besides $e_i$. By Proposition \ref{properties_of_mczs} (ii), $M_i$ contains a minimal conflicting subgraph. Since $r(M_i - e_i) = 0$, this minimal conflicting subgraph must contain $e_i$. This completes the proof.
\end{proof}
\end{proposition}

From Proposition \ref{r_distinct_mczs}, it is clear that the resistance of a cubic class two graph is less than or equal to the number of distinct minimal conflicting subgraphs contained in the graph. We may be inclined to think that the number of minimal conflicting subgraphs is in some way upper bounded by resistance, however this is not the case. The flower snarks and Loupekine snarks represent counter examples to this idea. Each of the graphs in these classes have resistance 2. However, the order of the graphs can be arbitrarily large. Furthermore, the number of possible single vertices which can be removed from the said graphs in order to leave behind a minimal conflicting subgraph is also arbitrarily large. Thus the number of minimal conflicting subgraphs is not bounded by resistance. 

%While no such upper bound exists, there does exist an essential relationship between resistance and minimal conflicting subgraphs. The following theorem sheds much light this relationship, and will equip us to formally present a characterisation of the resistance of a graph in terms of the graph's minimal conflicting subgraphs. First, we present a necessary definition. Let $G$ be a subcubic class two graph with minimal conflicting subgraphs $M_1,\dots , M_r$. From each $M_i$ we may select a representative edge $e_i \in M_i$ to form a set $R = \{ e_1, \dots, e_r\}$. Such a set $R$ we call a \textit{representative conflicting subset} of $G$. Note that since some $e_i$ may be contained in more than one minimal conflicting subgraph, there may exist representative conflicting subsets with varying order. 

While no such upper bound exists, there does exist an essential relationship between resistance and minimal conflicting subgraphs. This relationship, as is proven in the following theorem, provides much insight on possible conflicting edges and minimal colourings of snarks. First, we present an important definition.

\begin{definition}
	Let $G$ be a subcubic class two graph with minimal conflicting subgraphs $M_1,\dots,M_r$. A representative conflicting subset of $G$ is a set of distinct edges $R = {e_1, \dots ,e_s} \subset E(G)$ such that $R \cap M_i \neq \emptyset$ for each $i$.
\end{definition}

We note that for a graph $G$ there may exist representative conflicting subsets of varying order.

\begin{theorem} \label{choose_confledges_in_mincol}
	%Let $G$ be a subcubic class two graph and let $R$ be a representative conflicting subset of $G$. Then $G-R$ is 3-edge-colourable. Furthermore, if $|R|$ is minimal then $r(G)=|R|$.
	Let $G$ be a subcubic class two graph. Then 
	$$r(G) = \min\{|R| : R~ {\rm is~ a~ representative~ conflicting~ subset~ of}~ G\}$$
	\begin{proof}
		Let $\mathcal{M} = \{M_1,\dots,M_m\}$ be the set of all minimal conflicting subgraphs in $G$ and let $R$ be a representative conflicting subset of $G$. Note that no $M_i$ in $\mathcal{M}$ is a subgraph of $G-R$. Assume now that $G-R$ is not 3-edge-colourable. Then $G-R$ contains some minimal conflicting subgraph $M'$ by Proposition \ref{properties_of_mczs} (ii). But $M'$ is also contained in $G$, which is a contradiction since $M'$ is not contained in $\mathcal{M}$. Therefore $G-R$ is 3-edge-colourable.
		
		Let $|R|$ be minimal. Since $G-R$ is 3-edge-colourable, we know that $r(G) \leq |R|$. Assume that $r(G) < |R|$ and let $f$ be a minimal colouring of $G$. Let $R'$ be the conflicting edges with regard to $f$. By Proposition \ref{r_distinct_mczs}, every element in $R'$ is contained in some minimal coflicting subgraph of $G$. If $R'$ is not a representative conflicting subset then there exists some minimal conflicting subgraph $M' \subset G$ which contains no conflicting edges with regard to $f$. In which case, we have a minimal conflicting subgraph of $G$ which is properly coloured by $f$ using just three colours, a contradiction. If $R'$ is a representative conflicting subset, then the minimality of $|R|$ is contradicted since $|R'| = r(G) < |R|$. Therefore, $r(G)=|R|$.
	\end{proof}
\end{theorem}

%\begin{remark} \label{rem_characterising_resistance} {\rm
%Using Theorem \ref{choose_confledges_in_mincol}, we are now able to explicitly characterise the relationship between the resistance of a cubic graph and its  minimal conflicting subgraphs. The resistance of a graph is equal to the minimum number of non-empty intersections of minimal conflicting subgraphs such that the union of these intersections itself has nonempty intersection with each minimal conflicting subgraph. This is because we can select one edge from each of these nonempty intersections to form a representative conflicting subset of minimal order. }
%\end{remark}

%\begin{theorem} \label{thm_characterising_resistance} 
%Let $G$ be a subcubic class two graph and let $\mathcal{M} = \{M_1,\dots,M_m\}$ be the set of all minimal conflicting subgraphs in $G$. Let $\mathcal{I}$ be the set of all non-empty intersections of one or more elements of $\mathcal{M}$. Then   
%$$r(G) = \min \{ r ~| ~{\rm there~ exists}~ I_1,\dots,I_r \in \mathcal{I} ~{\rm with}~ (I_1 \cup \dots \cup I_r) \cap M_i \neq \varnothing {\rm~for~ every}~M_i \in \mathcal{M} \} .$$  
%\begin{proof}
%This result follows on directly from Theorem \ref{choose_confledges_in_mincol} and Remark \ref{rem_characterising_resistance}
%\end{proof}
%\end{theorem}

A cubic graph $G$ may have resistance $r(G)$, but given that information there is no way of knowing which combination of $r(G)$ edges may be removed from $G$ in order to render colourability. Theorem \ref{choose_confledges_in_mincol} is significant in that it informs us exactly which combinations of $r(G)$ edges are sufficient for this purpose. The requisite is that we can identify the minimal conflicting subgraphs of $G$. Another way of understanding the result, is that we can choose a minimal colouring, relative to conflicting edges, by simply selecting a combination of edges from each minimal conflicting subgraph, as long as this is done minimally.
		
%Theorem \ref{choose_confledges_in_mincol} is further significant in that it shows us that choosing a minimal colouring, relative to equivalence, is a simple matter of selecting a combination of edges from each minimal conflicting subgraph, as long as this is done minimally. 
Furthermore, with Theorem \ref{choose_confledges_in_mincol} we further note that if there exists some edge $e$ which is contained in exactly one minimal conflicting subgraph $M$, but $M$ has non-empty intersection with some other minimal conflicting subgraph $M'$, then $e$ may not be conflicting in any minimal colouring of $G$. Another way of saying this is, if $e \in M$ where $M$ is a minimal conflicting subgraph of $G$, then it is not necessarily the case that $r(G-e)=r(G)-1$. Equivalently, we could say that there does not necessarily exist some minimal colouring of $H \subset G$ which can be extended to a minimal colouring of $G$. It is possible to have a minimal colouring of a subgraph $H \subset G$ with say, $r_1$ conflicting edges, such that a minimal extension has $r_2$ further conflicting edges, but $r_1 + r_2 > r(G)$. What is also clear is that that if every minimal conflicting subgraph of $G$ is disjoint, then $r(G) = |\mathcal{M}|$, where $\mathcal{M}$ is the set of all minimal conflicting subgraphs in $G$. Consequent to this discussion, we define the following.

\begin{definition} \label{critical_zone} {\rm
		Let $G$ be a subcubic graph. Let 
$$K_G = \{e \in G ~|~ f(e)=0 {\rm ~for ~some~ minimal~ colouring~}f{\rm~ of}~ G\}.$$ 
We call $K_G$ the \textit{critical subgraph} of $G$.}
\end{definition}

As we did with resistance, we are also able to explicitly characterise the critical subgraph in terms of the minimal conflicting subgraphs.

\begin{theorem}
%Let $G$ be a subcubic class two graph and let $\mathcal{M} = \{M_1,\dots,M_m\}$ be the set of all minimal conflicting subgraphs in $G$. Let $\mathcal{I}$ be the set of all non-empty intersections of one or more elements of $\mathcal{M}$. Then   
%$$K_G = \bigcup \{  I_1 \cup \dots \cup I_{r(G)} ~|~ (I_1 \cup \dots \cup I_{r(G)}) \cap M_i \neq \varnothing {\rm~for~ every}~ M_i \in \mathcal{M}\}.$$
Let $G$ be a subcubic class two graph. Then 
$$K_G = \bigcup\{R : R~ {\rm is~ a~ representative~ conflicting~ subset~ of}~G~{\rm of~minimal~order} \}.$$

\begin{proof}
%For any $I_1 \cup \dots \cup I_{r(G)}$ such that $(I_1 \cup \dots \cup I_{r(G)}) \cap M_i \neq \varnothing$ for every $M_i \in \mathcal{M}$, it is clear that every edge in $I_1 \cup \dots \cup I_{r(G)}$ is conflicting in some minimal colouring of $G$, by Remark \ref{rem_characterising_resistance}. Therefore, the union of all such  $I_1 \cup \dots \cup I_{r(G)}$ is contained in $K_G$.

%Let $f$ be a minimal colouring of $G$. Let $R = \{e_1, \dots, e_{r(G)} \}$ be the set of conflicting edges with regard to $f$. Then $R$ must be a representative conflicting subset of minimal order, by Theorem \ref{choose_confledges_in_mincol}. Therefore, we can find $I_1 , \dots , I_{r(G)} \in \mathcal{I}$ such that $e_1 \in I_1, \dots, e_{r(G)} \in I_{r(G)}$, and $(I_1 \cup \dots \cup I_{r(G)}) \cap M_i \neq \varnothing$ for every $M_i \in \mathcal{M}$, by Theorem \ref{thm_characterising_resistance}. Therefore, $K_G$ is contained in the union of all such $I_1 \cup \dots \cup I_{r(G)}$. This completes the proof.

Let $R$ be a representative conflicting subset of $G$ with minimal order. Then the edges in $R$ are the conflicting edges of some minimal colouring of $G$. Therefore, the edges in $R$ are all critical.

Let $e$ be a critical edge of $G$. Then it is a conflicting edge in some minimal colouring $f$ of $G$. Let $R$ be the set of all conflicting edges in $G$ with regard to $f$. Since $G - R$ is colourable, $R$ is a representative conflicting subset of $G$. Since $f$ is minimal, $R$ must have minimal order. Therefore, $e$ is contained in the union of all representative conflicting subsets of $G$ with minimal order.

\end{proof}
\end{theorem}

It is clear that $K_G \subseteq M_G$. %We present a simple example where $K_G$ is a strict subgraph of $M_G$, and also an example where $K_G = M_G$.
We present an example where $K_G \subset M_G = G$, an example where $K_G = M_G \subset G$, and an example where $K_G = M_G = G$. The first two examples are specific graphs, while the third is the interesting general case of a hypo-Hamiltonian snark.

\begin{example}  \label{ex1} {\rm
	The subcubic graph $G$ depicted below consists of four identical minimal conflicting subgraphs, $M_1, M_2, M_3$ and $M_4$. Thus $M_G = G$. $M_1 \cap M_2$, $M_2 \cap M_3$ and $M_3 \cap M_4$ are
	represented by the thicker edges. We have $r(G) = 2$ and $K_G = (M_1 \cap M_2) \cup (M_3 \cap M_4) \subset M_G = G$. The sets of two edges, one each from $(M_1 \cap M_2)$ and $(M_3 \cap M_4)$, are the only representative conflicting subsets of minimal order. Thus, even though $K_G = (M_1 \cap M_2) \cup (M_3 \cap M_4)$ is itself 3-edge-colourable, any minimal colouring of $G$ must contain a conflicting edge in each of $(M_1 \cap M_2)$ and $(M_3 \cap M_4)$.

	
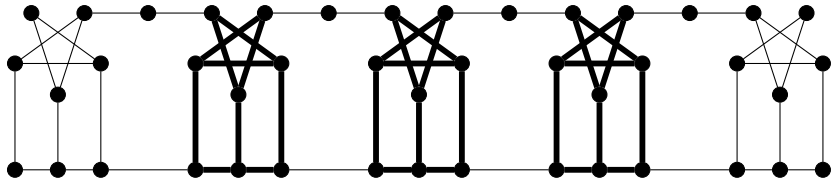
\begin{center}
		\begin{tikzpicture}[every node/.style={draw,shape=circle,fill=black,text=white,scale=0.6},scale=0.4]
		\begin{scope} %1
				\foreach \i [count=\ii from 0] in {54,126,198,270,342}{
					\path (\i:15mm) node (a\ii) {};
				}
				\foreach \x in {0,...,4}{
					\tikzmath{
						integer \y;
						\y = mod(\x+2,5);
					}
				\draw (a\x) -- (a\y);
				}
				\path let \p1 = (a2) in node (af0) at (\x1,-4) {};
				\path let \p1 = (a3) in node (af1) at (\x1,-4) {};
				\path let \p1 = (a4) in node (af2) at (\x1,-4) {};
				\draw (a2) -- (af0);
				\draw (a3) -- (af1);
				\draw (a4) -- (af2);
				\draw (af0) -- (af1) -- (af2);
			\end{scope}
			\begin{scope} [xshift=6cm]%2
				\foreach \i [count=\ii from 0] in {54,126,198,270,342}{
					\path (\i:15mm) node (b\ii) {};
				}
				\foreach \x in {0,...,4}{
					\tikzmath{
						integer \y;
						\y = mod(\x+2,5);
					}
				\draw (b\x) -- (b\y) [line width=2.2pt];
				}
				\path let \p1 = (b2) in node (bf0) at (\x1,-4) {};
				\path let \p1 = (b3) in node (bf1) at (\x1,-4) {};
				\path let \p1 = (b4) in node (bf2) at (\x1,-4) {};
				\draw (b2) -- (bf0)  [line width=2.2pt] ;
				\draw (b3) -- (bf1)  [line width=2.2pt];
				\draw (b4) -- (bf2) [line width=2.2pt];
				\draw (bf0) -- (bf1) -- (bf2) [line width=2.2pt];
			\end{scope}
			\begin{scope} [xshift=12cm]%3
				\foreach \i [count=\ii from 0] in {54,126,198,270,342}{
					\path (\i:15mm) node (c\ii) {};
				}
				\foreach \x in {0,...,4}{
					\tikzmath{
						integer \y;
						\y = mod(\x+2,5);
					}
				\draw (c\x) -- (c\y)  [line width=2.2pt];
				}
				\path let \p1 = (c2) in node (cf0) at (\x1,-4) {};
				\path let \p1 = (c3) in node (cf1) at (\x1,-4) {};
				\path let \p1 = (c4) in node (cf2) at (\x1,-4) {};
				\draw (c2) -- (cf0)  [line width=2.2pt];
				\draw (c3) -- (cf1)  [line width=2.2pt];
				\draw (c4) -- (cf2)  [line width=2.2pt];
				\draw (cf0) -- (cf1) -- (cf2)  [line width=2.2pt];
			\end{scope}
			\begin{scope} [xshift=18cm]%3
				\foreach \i [count=\ii from 0] in {54,126,198,270,342}{
					\path (\i:15mm) node (d\ii) {};
				}
				\foreach \x in {0,...,4}{
					\tikzmath{
						integer \y;
						\y = mod(\x+2,5);
					}
					\draw (d\x) -- (d\y)  [line width=2.2pt];
				}
				\path let \p1 = (d2) in node (df0) at (\x1,-4) {};
				\path let \p1 = (d3) in node (df1) at (\x1,-4) {};
				\path let \p1 = (d4) in node (df2) at (\x1,-4) {};
				\draw (d2) -- (df0)  [line width=2.2pt];
				\draw (d3) -- (df1)  [line width=2.2pt];
				\draw (d4) -- (df2)  [line width=2.2pt];
				\draw (df0) -- (df1) -- (df2)  [line width=2.2pt];
			\end{scope}
			\begin{scope} [xshift=24cm]%3
				\foreach \i [count=\ii from 0] in {54,126,198,270,342}{
					\path (\i:15mm) node (e\ii) {};
				}
				\foreach \x in {0,...,4}{
					\tikzmath{
						integer \y;
						\y = mod(\x+2,5);
					}
					\draw (e\x) -- (e\y);
				}
				\path let \p1 = (e2) in node (ef0) at (\x1,-4) {};
				\path let \p1 = (e3) in node (ef1) at (\x1,-4) {};
				\path let \p1 = (e4) in node (ef2) at (\x1,-4) {};
				\draw (e2) -- (ef0);
				\draw (e3) -- (ef1);
				\draw (e4) -- (ef2);
				\draw (ef0) -- (ef1) -- (ef2);
			\end{scope}
		
		\draw (a0) -- (b1) node (mid) [midway] {};
		\draw (b0) -- (c1) node (mid) [midway] {};
		\draw (c0) -- (d1) node (mid) [midway] {};
		\draw (d0) -- (e1) node (mid) [midway] {};
		\draw (af2) -- (bf0);
		\draw (bf2) -- (cf0);
		\draw (cf2) -- (df0);
		\draw (df2) -- (ef0);
			
		\end{tikzpicture}
	\end{center}}
\end{example}

\begin{example}  \label{ex2} {\rm
	The snark $G$ depicted below consists of three identical  non-overlapping minimal conflicting subgraphs $M_1, M_2$ and $M_3$. $M_1 \cup M_2 \cup M_3$ is represented by the thicker edges. Any set of three	edges, one each from $M_1, M_2$ and $M_3$, is a representative conflicting subset. Therefore $r(G) = 3$ and $K_G = M_1 \cup M_2 \cup M_3 = M_G \subset G.$
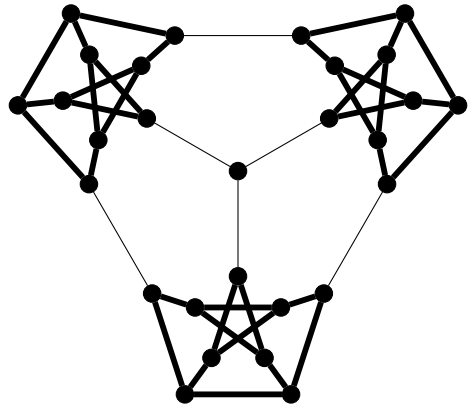
\begin{center}
	\begin{tikzpicture}[every node/.style={draw,shape=circle,fill=black,text=white,scale=0.7},scale=0.4]
		\begin{scope} [rotate=60, yshift=5cm, rotate=288]
			\foreach \i [count=\ii from 0] in {54,126,198,270,342}{
				\path (\i:15mm) node (a\ii) {};
			}
			\foreach \x in {0,...,4}{
				\tikzmath{
					integer \y;
					\y = mod(\x+2,5);
				}
				\draw (a\x) -- (a\y) [line width=2.2pt];
			}
			\foreach \i [count=\ii from 0] in {54,126,198,270}{
				\path (\i:30mm) node (b\ii) {};
			}
			\draw (b0) -- (b1) [line width=2.2pt];
			\draw (b1) -- (b2) [line width=2.2pt];
			\draw (b2) -- (b3) [line width=2.2pt];
			\draw (a0) -- (b0) [line width=2.2pt];
			\draw (a1) -- (b1) [line width=2.2pt];
			\draw (a2) -- (b2) [line width=2.2pt];
			\draw (a3) -- (b3) [line width=2.2pt];
		\end{scope} 
		\begin{scope} [rotate=180, yshift=5cm, rotate=288] % [xshift=8cm, rotate=234]%1 
			\foreach \i [count=\ii from 0] in {54,126,198,270,342}{
				\path (\i:15mm) node (c\ii) {};
			}
			\foreach \x in {0,...,4}{
				\tikzmath{
					integer \y;
					\y = mod(\x+2,5);
				}
				\draw (c\x) -- (c\y) [line width=2.2pt];
			}
			\foreach \i [count=\ii from 0] in {54,126,198,270}{
				\path (\i:30mm) node (d\ii) {};
			}
			\draw (d0) -- (d1) [line width=2.2pt];
			\draw (d1) -- (d2) [line width=2.2pt];
			\draw (d2) -- (d3) [line width=2.2pt];
			\draw (c0) -- (d0) [line width=2.2pt];
			\draw (c1) -- (d1) [line width=2.2pt];
			\draw (c2) -- (d2) [line width=2.2pt];
			\draw (c3) -- (d3) [line width=2.2pt];
		\end{scope} 
		\begin{scope} [rotate=300, yshift=5cm, rotate=288]%[xshift=4cm, yshift=-6cm, rotate=108]%1 
			\foreach \i [count=\ii from 0] in {54,126,198,270,342}{
				\path (\i:15mm) node (e\ii) {};
			}
			\foreach \x in {0,...,4}{
				\tikzmath{
					integer \y;
					\y = mod(\x+2,5);
				}
				\draw (e\x) -- (e\y) [line width=2.2pt];
			}
			\foreach \i [count=\ii from 0] in {54,126,198,270}{
				\path (\i:30mm) node (f\ii) {};
			}
			\draw (f0) -- (f1) [line width=2.2pt];
			\draw (f1) -- (f2) [line width=2.2pt];
			\draw (f2) -- (f3) [line width=2.2pt];
			\draw (e0) -- (f0) [line width=2.2pt];
			\draw (e1) -- (f1) [line width=2.2pt];
			\draw (e2) -- (f2) [line width=2.2pt];
			\draw (e3) -- (f3) [line width=2.2pt];
		\end{scope} 
		\node (c) at (0,0) {};
		\draw (a4) -- (c);
		\draw (c4) -- (c);
		\draw (e4) -- (c);
		\draw (b3) -- (d0);
		\draw (d3) -- (f0);
		\draw (f3) -- (b0);
	\end{tikzpicture}
	\end{center}}
\end{example}

\begin{example} \label{ex3} {\rm
	The general case of a hypo-Hamiltonian snark $G$ is depicted below. Let $e$ be any given edge in $G$. Then $e$ is contained in a Hamiltonian cycle of $G-v$ where $v$ is a vertex distance 1 from $e$. In the diagram, $e$ is conflicting in a minimal 4-edge-colouring of $G$ with two conflicting edges. Therefore, $K_G = G$. The 3-coloured chordal edges and the alternatively 1-2 coloured edges in the Hamiltonian cycle are not depicted in the diagram. Since hypo-Hamiltonian snarks are bicritical \cite{steffen2}, we note that this implies that every minimal conflicting subgraph of a hypo-Hamiltonian snark contains all but one vertex.
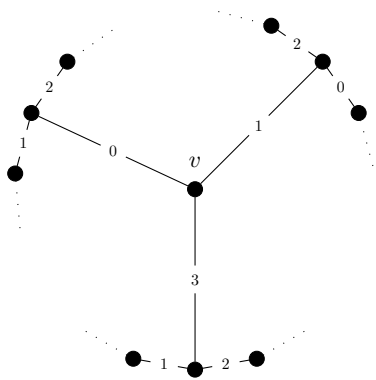
\begin{center}
\begin{tikzpicture}[every node/.style={draw,shape=circle,fill=black,text=white,scale=0.6},scale=0.6]

				\node (a0) at (25:4cm) [label={[above left,color=black,scale=1.5]:$ $}] {};
				\node (a1) at (45:4cm) [label={[above right,color=black,scale=1.5]:$ $}] {};
				\node (a2) at (65:4cm) [label={[above left,color=black,scale=1.5]:$ $}] {};

				\node (a3) at (135:4cm) [label={[above left,color=black,scale=1.5]:$ $}] {};
				\node (a4) at (155:4cm) [label={[above left,color=black,scale=1.5]:$ $}] {};
				\node (a5) at (175:4cm) [label={[above left,color=black,scale=1.5]:$ $}] {};
				
				\node (a6) at (250:4cm) [label={[above left,color=black,scale=1.5]:$ $}] {};
				\node (a7) at (270:4cm) [label={[above left,color=black,scale=1.5]:$ $}] {};
				\node (a8) at (290:4cm) [label={[above left,color=black,scale=1.5]:$ $}] {};

				\draw node (x1) at  (5:4cm) [fill=white,text=black,draw = none,scale=1] {} ;
				\draw node (x2) at  (85:4cm) [fill=white,text=black,draw = none,scale=1] {} ;
				\draw node (x3) at  (115:4cm) [fill=white,text=black,draw = none,scale=1] {} ;
				\draw node (x4) at  (195:4cm) [fill=white,text=black,draw = none,scale=1] {} ;
				\draw node (x5) at  (230:4cm) [fill=white,text=black,draw = none,scale=1] {} ;
				\draw node (x6) at  (310:4cm) [fill=white,text=black,draw = none,scale=1] {} ;
				
				\node (c) at (0,0) [label={[above,color=black,scale=1.5]:$v$}] {};
				
				\draw (c) -- node [draw=none,fill=white,text=black] {$1$} (a1);
				\draw (c) -- node [draw=none,fill=white,text=black] {$3$} (a7);
				\draw (x1) -- (a0) [loosely dotted];
				\draw (x2) -- (a2) [loosely dotted];
				\draw (x3) -- (a3) [loosely dotted];
				\draw (x4) -- (a5) [loosely dotted];
				\draw (a4) -- node [draw=none,fill=white,text=black] {$0$} (c);
				\draw (x5) -- (a6) [loosely dotted];
				\draw (x6) -- (a8) [loosely dotted];
				\draw (a0) -- node [draw=none,fill=white,text=black] {$0$} (a1) -- node [draw=none,fill=white,text=black] {$2$}(a2);
				\draw (a3) -- node [draw=none,fill=white,text=black] {$2$} (a4) -- node [draw=none,fill=white,text=black] {$1$} (a5);
				\draw (a6) -- node [draw=none,fill=white,text=black] {$1$} (a7) -- node [draw=none,fill=white,text=black] {$2$} (a8);

			\end{tikzpicture}
\end{center}}
\end{example} 

%%%%%%%%%%%%%%%%%%%%%%%%%%%%%%%%%%%%%%%%%%%%%%%%%%%%%%%%%%%%%%
%%%%%%%%%%%%%%%%%%%%%%%%%%%%%%%%%%%%%%%%%%%%%%%%%%%%%%%%%%%%%%
%%%%%%%%%%%%%%%%%%%%%%%%%%%%%%%%%%%%%%%%%%%%%%%%%%%%%%%%%%%%%%
%%%%%%%%%%%%%%%%%%%%%%%%%%%%%%%%%%%%%%%%%%%%%%%%%%%%%%%%%%%%%%
%%%%%%%%%%%%%%%%%%%%%%%%%%%%%%%%%%%%%%%%%%%%%%%%%%%%%%%%%%%%%%
%%%%%%%%%%%%%%%%%%%%%%%%%%%%%%%%%%%%%%%%%%%%%%%%%%%%%%%%%%%%%%

\section{Further considerations}

\subsection{Clusters}

We have noticed that it is typically the case in smaller snarks that minimal conflicting subgraphs have non-empty intersections. To facilitate further brief discussion, it serves to formally define groups of minimal conflicting subgraphs in terms of non-empty intersections, as well as distinguish between different types of these groups.

\begin{definition} \label{cluster} {\rm 
		Let $G$ be a subcubic class two graph. Let $\mathcal{M} = \{ M_1, \dots, M_m \}$ be a collection of minimal conflicting subgraphs of $G$. 
\begin{enumerate}
\item [(i)] If for every $i \in \{1,\dots,m\}$ with $i\neq j$ there exists some $j \in \{1,\dots,m\}$ such that $M_i \cap M_j \neq \varnothing$, and $M \cap M_i = \varnothing$ for any other minimal conflicting subgraph $M \notin \mathcal{M}$, then we call $\mathcal{M}$ a \textit{cluster} of minimal conflicting subgraphs. 
\item [(ii)] If $\mathcal{M}$ is a cluster and $\bigcap M_i \neq \varnothing$ then we call $\mathcal{M}$ a \textit{dense cluster}.
\item [(iii)] If $\mathcal{M}$ is a cluster and is not dense then it is a \textit{sparse cluster}. 
\item [(iv)]  If $\mathcal{M}$ is sparse cluster such that for every $i, j \in \{1,\dots,m\}$ we have that $M_i \cap M_j \neq \varnothing$, then it is a \textit{densely sparse cluster}.
\end{enumerate}}
\end{definition}

We prove and discuss some immediate results on these structures. Our investigations suggest that it serves to consider strictly subcubic clusters and cubic clusters seperately.

%\begin{proposition} \label{k(G)=1_no_standalone_mcz}
%	Let $G$ be a bridgeless cubic graph. If $M_G=G$ then $G$ consists entirely of one sparse cluster of minimal conflicting subgraphs.
%	\begin{proof}
%		Since the distance between any two clusters must be at least one by Proposition \ref{properties_of_mczs}, that one edge cannot be contained in any minimal conflicting subgraph. Thus $M_G=G$ implies that $G$ consists entirely of one cluster of minimal conflicting subgraphs. Assume the cluster is dense. Then there exists some $e \in \bigcap M_i$ so that $G -\{e\}$ contains no conflicting subgraphs. However, $G$ is then a cubic graph with resistance one, which is impossible. Therefore, $G$ consists entirely of one sparse cluster of minimal conflicting subgraphs.
%	\end{proof} 
%\end{proposition}	 

\begin{proposition} \label{propclusters}
The following statements are true.
\begin{enumerate}
	\item[(i)] There exists no cubic dense cluster.
	%\item[(ii)] There exists a strictly subcubic densely sparse cluster with 3 minimal conflicting subgraphs.	
	\item[(ii)] Let G be a bridgeless cubic graph. If $M_G = G$ then $G$ consists entirely of one sparse cluster of minimal conflicting zones.
	\item[(iii)] There exists a strictly subcubic %sparse 
	cluster with $n$ minimal conflicting subgraphs for each $n \geq 1$.
	%\item[(iv)] There exists a strictly subcubic sparse cluster with $n$ minimal conflicting subgraphs for each $n \geq 3$.
	%\item[(iv)] There exists cubic densely sparse clusters with more than 3 minimal conflicting subgraphs.

\end{enumerate}
\begin{proof} \hfill
	\begin{enumerate}
		\item[(i)] Every dense cluster has a representative conflicting subset of order 1. No cubic graph can have resistance 1. Therefore, no dense cluster can be cubic.
		\item[(v)] Since the distance between any two clusters must be at least one by Proposition \ref{properties_of_mczs}, that one edge cannot be contained in any minimal conflicting subgraph. Thus $M_G = G$ implies that $G$ consists entirely of one cluster of minimal conflicting subgraphs. By $(i)$, $G$ cannot be dense, and is therefore sparse.
		%\item[(ii)] Consider Example \ref{ex1}, but with three minimal conflicting subgraphs, and let $M1$ and $M3$ intersect in the same manner as $M1$ and $M2$, and $M2 and M3$. The result is a strictly subcubic densely sparse cluster with 3 minimal conflicting subgraphs.
		\item[(iii)] Consider Example \ref{ex1}, but with %$n \geq 4$ 
		$n$ minimal conflicting subgraphs $M_1, \dots , M_n$. As in Example \ref{ex1}, let $M_i$ intersect with $M_{i+1}$ for $i \in \{1, \dots , n-1\}$. The result is a strictly subcubic %sparse 
		cluster with $n$ minimal conflicting subgraphs.
		%\item[(iv)] As in (iii), the cases where $n \geq 3$ is sparse.
		%\item[(iv)] Any hypo-Hamiltonian graph, as in Example \ref{ex3}, is a cubic densely sparse clusters with more than 3 minimal conflicting subgraphs.
	\end{enumerate}
	\end{proof}
\end{proposition}

%The nature of sparse clusters is intriguing. In particular, the difference between sparse clusters and densely sparse clusters. Given our investigations and Proposition \ref{propclusters}, we suspect that if a sparse cluster with more than three minimal conflicting subgraphs is densely sparse, then it is cubic. 

Now, any cluster with two minimal conflicting subgraphs is trivially dense, and we have seen that we can easily find such clusters. In Proposition \ref{propclusters} (iii), the clusters are however sparse for $n \geq 3$. The question of whether there exists dense clusters with three or more minimal conflicting subgraphs remains.  

\begin{problem} \label{prob2}
	For which $n \geq 3$ does there exist a dense cluster with $n$ minimal conflicting subgraphs?
\end{problem}

%\begin{problem} \label{prob1}
%	Does there exist cubic sparse clusters with 3 minimal conflicting subgraphs?
%\end{problem}

Recall Example \ref{ex3}, that in a hypo-Hamiltonian snark, the removal of any vertex leaves behind a minimal conflicting subgraph containing all the remaining vertices. Thus any two minimal conflicting subgraphs intersect, and there is no single vertex which is present in every minimal conflicting subgraph. Hypo-Hamiltonian snarks are therefore densely sparse clusters. From our investigations, we suspect that the only cubic densely sparse clusters are those which are similar to Example \ref{ex3}. That is, possibly hypo-Hamiltonian, and consequently where the removal of any one vertex leaves behind a minimal conflicting subgraph. In such cases as well, resistance is necessarily $2$. Recall as well that every edge in a ypo-Hamiltonian snark is critical.
%If a sparse cluster is strictly subcubic, then it possibly allows for the addition of edges without adding more minimal conflicting subgraphs so that $M_G$ is strictly contained in $G$. Thus in the hypothesis of Proposition \ref{propclusters} (v), we suspect that the sparse cluster is necessarily densely sparse. 
Thus, we formulate the following conjecture.

\begin{conjecture} \label{con2}
	Let $G$ be a bridgeless cubic graph. Then $K_G = G$ if and only if $r(G) = 2$ and $G$ is a densely sparse cluster.
\end{conjecture}

%Obviously, Conjecture \ref{m=2 implies r=2} implies Conjecture \ref{k=2 implies r=2}.
\subsection{Snark reduction}

	Although triviality in snarks is not formally defined, snarks have generally been considered to contain more triviality if they are easily reducible to smaller snarks by some well-defined reduction (many variations of snark reductions have been considered by previous authors, (see for example \cite{nedela, steffen2})). In particular, the most universally accepted notion of triviality in snarks are those with either girth less than 5 or those with cyclic connectivity less than 4. 
	There does however exist many snarks with these ``trivial" properties, with arbitrarily large resistance, and which can only be reduced to smaller snarks with strictly less resistance. Reduction of resistance cannot be considered trivial, since some structural complexity is contributing to that resistance being present. Thus we propose instead that snarks which can be reduced in size, without reducing resistance, should be considered more trivial, and that so-called ``trivial'' snarks should be studied as much as so-called ``non-trivial'' snarks. Given Theorem \ref{choose_confledges_in_mincol}, these would typically be snarks with a large buffer subgraph. Thus we may perhaps formally define a snark to contain no triviality if its buffer subgraph is empty.
	%Snarks with girth less than 5 and cyclic connectivity less than 4 have generally been considered to be more trivial than other snarks. That is because these snarks are easily reducible to smaller snarks. There are however many snarks without these characteristics which contain large buffer subgraphs, and as such should be considered to contain much triviality. We propose instead that snarks which can be reduced in size, without reducing resistance, should be considered more trivial. 	These are typically cases of snarks with large buffer subgraphs. 
	The notion of a maximal conflicting subgraph and buffer subgraph therefore opens up a new avenue of consideration regarding snark reductions. That is, reducing the snark to contain only the essentially uncolourable, by removing the buffer subgraph.
	
	This new notion of reducibility relates interestingly to a problem of oddness and resistance in snarks. Recall that the oddness of a graph $G$ is the minimum number of odd components in a 2-factor of $G$, denoted as $\omega(G)$ \cite{huck}. In \cite{allie} we disproved a conjecture by Fiol et. al. in \cite{fiol}, by showing that the ratio of oddness to resistance can be arbitrarily large. It was conjectured in \cite{fiol} that $\omega(G) \leq 2r(G)$ for any snark $G$. We disproved this by constructing a class of graphs with increasing oddness, but constant resistance equal to 3. Interestingly, each of the graphs in the class defined contains exactly three disjoint minimal conflicting subgraphs (thus resistance is 3 in each graph). The increase in oddness, whilst keeping resistance constant, is as a result of adding particular subgraphs to the buffer subgraph. This leads to an interesting reformulation of the disproved conjecture which was posed in \cite{fiol}. 

\begin{conjecture} {\rm 
Let $G$ be a snark with an empty buffer subgraph. Then $\omega(G) \leq 2r(G)$.} 
\end{conjecture}


\begin{thebibliography}{99}                             

\bibitem{allie}
I.~Allie, \textit{Oddness to resistance ratios in cubic graphs}, Discrete Math. 342 (2019), pp. 387 - 391.

\bibitem{bondy}
J. A. Bondy, U. S. R. Murty, \textit{Graph Theory with Applications},  Elsevier Science Publishing Co., New York  (1976). 

\bibitem{fiol} 
M. A. Fiol, G. Mazzuoccolo, E. Steffen, \textit{Measures of Edge-Uncolorability of Cubic Graphs}, Electron. J. Combin. 25(4) (2018), P4.54. 

\bibitem{huck}
A. Huck, M. Kochol, \textit{Five cycle double covers of some cubic graphs}, J. Combin. Theory Ser. B 64 (1995), pp. 119 - 125. 

\bibitem{lukotka} 
R. Lukot'ka, J. Maz\'{a}k, \textit{Weak oddness as an approximation off oddness and resistance in cubic graphs}, (2016), arXiv:1602.02949.

\bibitem{nedela}
R. Nedela, M. Skoviera, \textit{Decompositions and reductions of snarks}, Journal of Graph Theory
22 (1996), pp. 253–279.

\bibitem{steffen} 
E. Steffen, \textit{Measurements of edge-uncolorability}, Discrete Math. 280 (2004), pp. 191 - 214.

\bibitem{steffen2} 
E. Steffen, \textit{Classifications and characterizations of snarks}, Discrete Math. 188 (1998), pp. 183 - 203.

\bibitem{szekeres}
G. Szekeres,  \textit{Polyhedral decomposition of cubic graphs}, Bull. Aust. Math. Soc. 8 (1973), pp. 367-387.

\bibitem{tutte}
W. T. Tutte,  \textit{A contribution to the theory of chromatic polynomials}, Canadian J. Math. 6 (1954), pp. 80-91.


\end{thebibliography}
\end{document}